\documentclass[a4paper]{article}
\usepackage{amsmath, amssymb,mathtools}
\usepackage{stmaryrd}
\usepackage{amsthm,cite}
\usepackage{amsfonts}
\usepackage{graphicx,subfigure}
\usepackage{epsfig,color,pifont,enumerate}
\usepackage{setspace}
\graphicspath{{img/}}

\usepackage{geometry}
\newcommand{\real}{\mathbb{R}}

\renewcommand{\phi}{\varphi}
\renewcommand{\epsilon}{\varepsilon}
\renewcommand{\kappa}{\varkappa}
\newcommand{\sgn}{\text{\bf sgn}\,}
\newcommand{\trs}{\top}
\newcommand{\br}{\mathbb{R}}
 \newcommand{\bldr}{\boldsymbol{r}}

\newtheorem{theorem}{Theorem}[section]

\newtheorem{lemma}{Lemma}[section]
\newtheorem{cor}{Corollary}[section]
\newtheorem{proposition}{Proposition}[section]

\newcommand{\so}{\text{\scriptsize $\mathcal{O}$}}
\newcommand{\spr}[2]{\left\langle #1; #2 \right\rangle}

\newcommand{\ov}[1]{\overline{#1}}

\newcommand{\dd}[1]{#1^{\prime\prime}}

\newcommand{\pf}{{\bf Proof:}\;}
\newcommand{\epf}{$\qquad \Box$}

\newcommand{\pp}{{\prime\prime}}
\title{Technical Facts About Dynamic Scalar Fields Underlying Algorithms of Mobile Robots Navigation for Tracking Environmental Boundaries and Extremum Seeking}
\date{}
\author{Alexey S. Matveev, Anna A Semakova, and Andrey V.Savkin}
\begin{document}
\maketitle
\section{Introduction}
Mobile robots
 are increasingly used for various missions in hazardous and complex environments not only to avoid or mitigate risks to human life and health but also due to their
lightweights, inexpensive components, and low power consumptions. These robots are often subjected to limitations on communication and so should be equipped with navigation and control systems that enable them to autonomously operate either constantly or for extended periods of time and distance.
Recently, such problems have gained much interest in control community. Tracking of environmental level sets is among them: The robot should reach and then track and so display the isoline where an unknown scalar field assumes a certain value. By doing so, the robot localizes the area inside this isoline where the field assumes greater values and which is typically of major interest. For example, this is a critically contaminated area if the field gives the radiation level or a pollutant concentration.
Such tasks are also of interest for tracking of oil and chemical spills \cite{ClaFie07}, detection and monitoring of harmful algae blooms \cite{PeDuJoPo12}, tracking forest fires \cite{CKBMLM06}
or contaminant clouds \cite{WTASZ05}, exploration of sea temperature, etc.
\par
Another mission of essential interest is extremum seeking, where a mobile robot should be autonomously navigated to the maximizer of an environmental field. In many typical scenarios, the field profile is unknown a priori and should be explored on-line based on point-wise measurements of the field. This holds in a variety of missions, such as environmental studies, tracking a moving target based on a single signal decaying away from the target, seeking sources of radiation emission or pollutant leakage, etc.
The discussed problem also falls into the general framework of extremum seeking control, which should ensure that a dynamic plant is directed to and subsequently tracks an optimal operating point based on on-line evaluation of the current performance, whose profile is unknown a priori  \cite{TaMoMaNeMa10,DoPeGu11}.
\par
Access to the field gradient is assumed by many works on robotic tracking of environmental boundaries (see e.g., \cite{MaBe03,BertKeMar04,SrRaKu08,ZhLe10,HsLoKu07} and literature therein) and environmental extremum seeking (see e.g., \cite{BaChSp00,TePo01,LiDu04,AtNyMiPa12,ZhOr12} and the literature therein).
However, spatial derivatives often are not measurable, unlike the field value at the robot's location. Their estimation needs simultaneous access to this value at several nearby locations distributed over two dimensions.
This may require complicated and costly hardware, e.g., several communicating robots and their ineffective (tight-flocking) formation, or many distant sensors on a single robot. Meanwhile, data exchange degradation due to e.g., communication constraints may require each robot to operate autonomously over considerable time and space even in multi-robot scenario, whereas equipment of a single robot with many distant sensors may constitute a separate problem, e.g., for small-size robots.
Anyhow, gradient estimators suffer from noise amplification, and their sensor-based implementation is still an intricate problem in practical setting \cite{AhAb07,chartrand11,VaKh08}.
This carries a threat of performance degradation, puts extra burden on controller tuning, and may greatly increase the overall computational load \cite{AhAb07}.
\par
Motivated by the foregoing, a substantial research activity was devoted to gradient-free approaches and scenarios where a single mobile sensor has a pointwise access to the field value only when tracking an environmental boundary (see e.g., \cite{BaRe03,CKBMLM06,Anders07,ZhBer07,JAHB09,CLBMM05,BarBail07,MaTeSa12,MaChaSa12a,MaHoOvAnSa15} and literature therein) or when looking for the maximizer of an environmental field (see e.g., \cite{BuYoBrSi96,ZhArGhSiKr07,CoKr09,CoSiGhKr09,LiKr10,ElBr12,BaLe02,MaSaTe08,BarBail08,MaTeSa11,MaSaTe08,MaTeSa11,GaZh03,MaHoSa3d14,MaHoSa15} and literature therein). These works mostly offer reactive controllers which directly convert current observation into current control input via a simple, reflex-like rule (for recent extensive surveys of reactive algorithms for navigation of mobile robots, we refer the reader to \cite{HoMaSa15rob,MaSaHoWa16b}).
\par
A common feature of the previous theoretical developments concerned with gradient-free tracking of environmental boundaries and extremum seeking navigation is that they deal mainly with steady fields. However in the real world, environmental fields are almost never steady and often cannot be well approximated by steady fields. The problem of extremum seeking in dynamic fields is also concerned with localization and circumnavigation of an unknowingly maneuvering target based on a single measurement that decays as the sensor goes away from the target, like the strength of the infrared, acoustic, or electromagnetic signal, or minus the distance to the target. Such navigation is of interest in many areas \cite{ADB04,GS04,MTS11ronly}; it carries a potential to reduce the hardware complexity and cost and improve target pursuit reliability. A solution for such problem in the very special case of the unsteady field --- minus the distance to an unknownly moving Dubins-car-like target --- was proposed and justified in \cite{MTS11ronly}. However the results of \cite{MTS11ronly} are not applicable to more general dynamic fields. To the best of the knowledge of the authors, rigorously justified navigation laws for gradient-free tracking of environmental boundaries of generic unsteady fields and for gradient-free extremum seeking in such fields are presented in \cite{MaHoOvAnSa15} and \cite{MaHoSa15}, respectively. However, \cite{MaHoOvAnSa15} follows the sliding mode paradigm, which gives rise to concerns about chattering and excessive control activity, and also employs the time-derivative of the field reading, whereas derivative estimators are prone to noise amplification. On the other hand, \cite{MaHoSa15} does not take into account possible nonholonomy of the robot and other kinematic or dynamic constraints that may apply.
\par
Thus in research concerned with unsteady environmental fields, many issues remain open. Their elaboration
inevitably demands for evaluation of the proposed navigation solutions versus the field properties.
In turn, this calls for systematic employment of relevant  characteristics of unsteady fields and their isolines that properly inform about their geometry, kinematics, deformations, etc., and also enjoy the benefit of being physically meaningful. General geometric and kinematic relations
also belong to expectedly needed armamentarium.
Meanwhile among the survey of the relevant literature, the authors failed to come across a systematic and comprehensive exposition of these issues.
\par
This paper is aimed at filling this gap. It presents some field characteristics and their properties that are relevant to research on environmental extremum seeking and boundary tracking, according to the experience of the authors. The proofs are mostly via rather boring, though straightforward calculus computations. This partly discloses the idea behind this text: To unload research papers on the above main topics from the need to perform these technical developments. A first example is the paper submitted by the authors to the IFAC journal {\it Automatica}, which refers to
the technical facts from this text.
\par
This paper is organized as follows. Section~\ref{sec1} gives a quick idea about the ultimate target of the paper, as is viewed by the authors: It sets up the problems of tracking environmental boundaries and environmental extremum seeking. Section~\ref{subsec.not} introduces and discusses some characteristics of dynamic scalar fields; Sect.~\ref{robot.deriv} studies time derivatives of the field reading of a moving robot.  
The last Section~\ref{init.sec.state} offers a useful technical estimation of deviation of a perpetually rotating robot from the initial state. 
\par
The following notations are used throughout the paper:
 \begin{itemize}
 \item the symbol $:=$ means ``is defined to be'';
  \item $\|\cdot\|$ is the Euclidian norm;
   \item $\spr{\cdot}{\cdot}$ is the standard inner product in the plane;
    \item $^\top$ stands for transposition;
     \item $\overset{\text{``A''}}{=}$ means that the equation $=$ is true by virtue of ``A'', similar convention is adopted for other binary relations;
        \item $R_{\alpha}$ is the counter clockwise rotation through angle $\alpha$.
\end{itemize}

\section{The Problems of Gradient-Free Isoline Tracking and Extremum Seeking by a Mobile Robot}\label{sec1}
We consider a point-wise robot traveling in a Cartesian plane with the absolute
coordinates $x$ and $y$. For the sake of definiteness and simplicity, we assume the simple integrator model of the robot's kinematics:
The robot is controlled by the time-varying linear velocity $\vec{v}$ whose magnitude does not exceed a given constant $\ov{v}$. The plane
hosts an unknown scalar time-varying field $D(t,\boldsymbol{r}) \in \br$, where $t$ is time and
$\boldsymbol{r}:= (x,y)^\trs$.
The on-board control system has access only to the field
value $d(t):= D[t,\bldr(t)]$ at the robot's current location $\bldr(t) = [x(t),y(t)]^\trs$.
No data about the derivatives of $D$ are available; in particular, the robot is not aware of the field's spatial gradient.
\par
The kinematic model of the robot is as follows:
\begin{equation}
\label{1}
\dot{\bldr} = \vec{v}, \qquad \bldr (0) = \bldr_{\text{in}}\qquad \|\vec{v}\| \leq \ov{v},
\end{equation}
where $\vec{v}$ is the control input. By introducing the speed $v := \|\vec{v}\|$ and velocity orientation angle $\theta$, this can be shaped into
\begin{equation}
\label{1aa}
\dot{\bldr} = v \vec{e}(\theta), \quad 0 \leq v \leq \ov{v}, \quad \text{where} \quad \vec{e}(\theta) := (\cos \theta, \sin \theta)^\trs
\end{equation}
and $v, \theta$ are control inputs.
\par
\begin{figure}
\centering
\subfigure[]{\scalebox{0.3}{\includegraphics{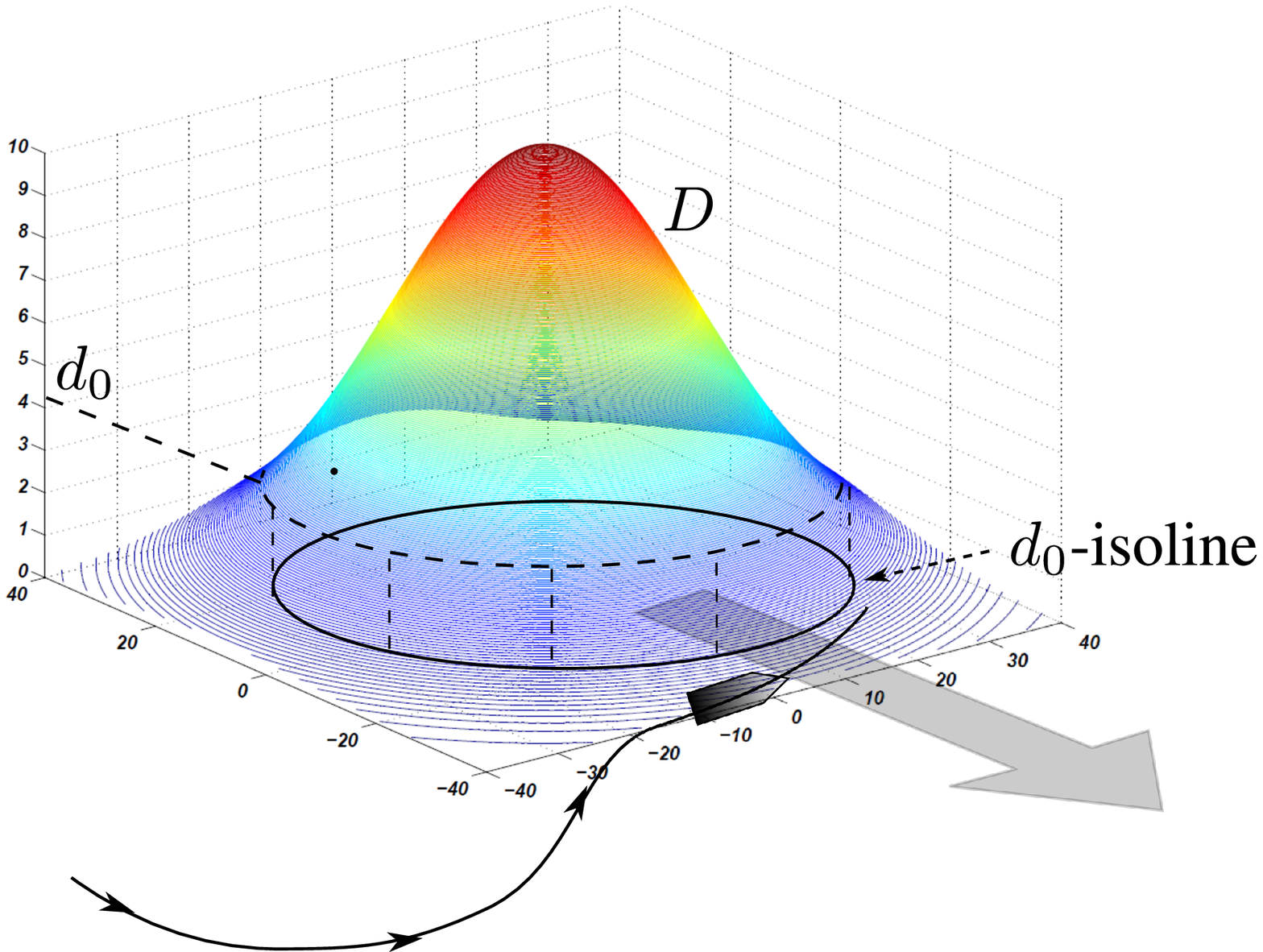}}}
\makebox[2.0cm][c]{}
\subfigure[]{\scalebox{0.3}{\includegraphics{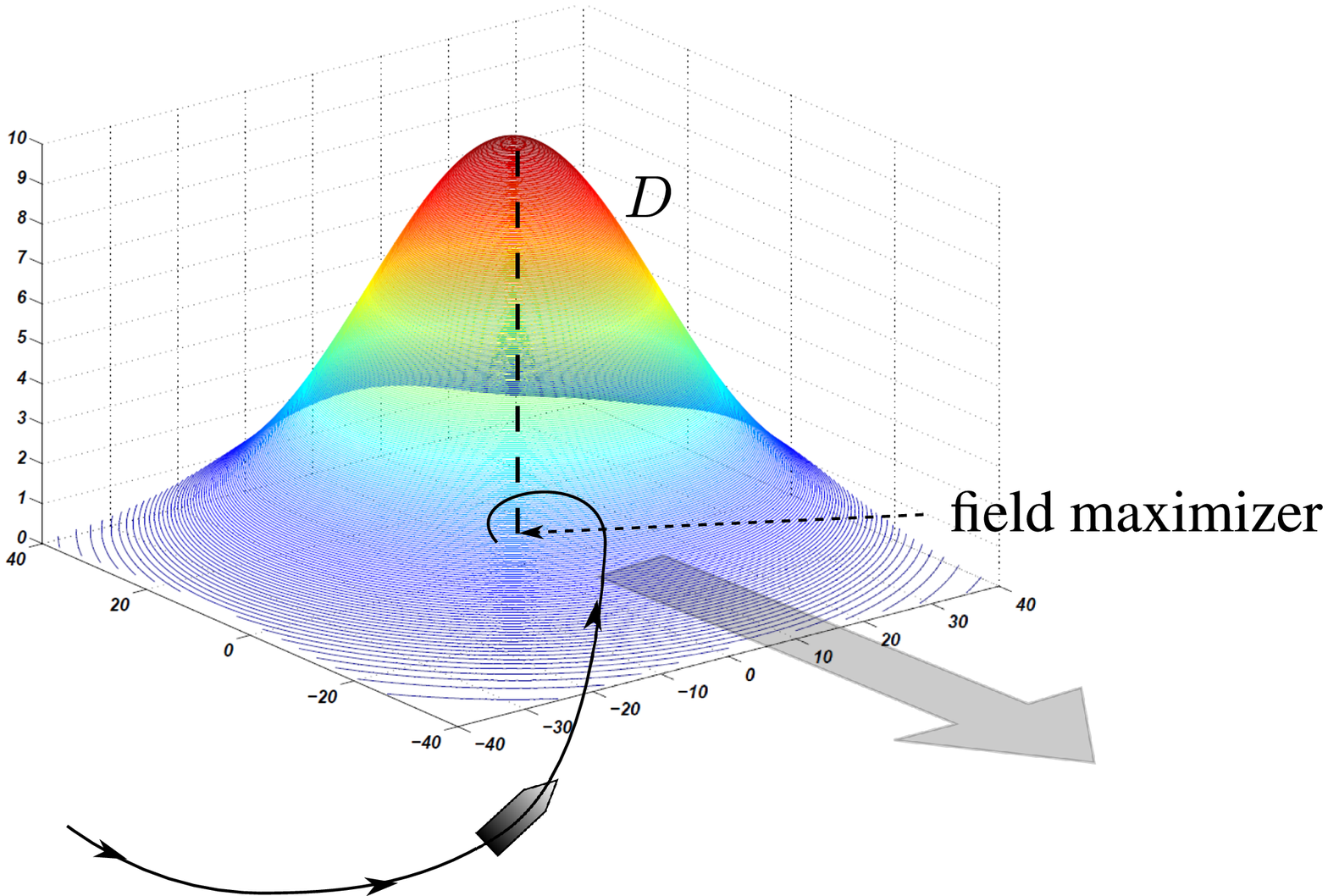}}}
\caption{(a) Tracking of $d_0$-isoline of a dynamic scalar field $D$; (b) Extremum seeking by using a mobile robot.}
\label{isol.fig}
\end{figure}

{\bf Isoline tracking problem.}
The robot should be steered to the level curve $D(t,\bldr) = d_0$ where the
field assumes a pre-specified value $d_0$ and to subsequently circulate along this curve; see Fig.~\ref{isol.fig}(a).
\par
{\bf Extremum seeking problem.} The objective is to drive
the robot to the moving point $\bldr^0(t)$ where
$D(t,\bldr)$ attains its maximum over $\bldr \in \br^2$ and then to keep the robot in a vicinity of
$\bldr^0(t)$, thus displaying the approximate location of
$\bldr^0(t)$; see Fig.~\ref{isol.fig}(b).
\section{Some Characteristics of a Dynamic Scalar Field}
\label{subsec.not}
\setcounter{equation}{0}
\begin{figure}[h]
\centering
\scalebox{0.3}{\includegraphics{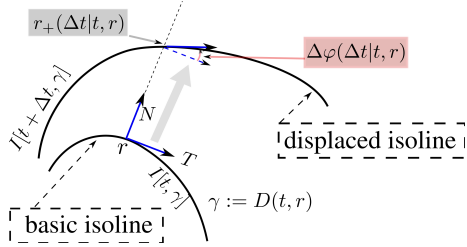}}
\caption{Two close isolines.
}
\label{fig.converge}
\end{figure}
We start with introducing notations and examined characteristics of dynamic fields.
\begin{itemize}
\item $\nabla = \left( \frac{\partial}{\partial x} , \frac{\partial}{\partial y} \right)^\top$ --- the spatial gradient;
    \item $D^{\prime\prime}$ --- the spatial Hessian, i.e., the matrix of the second derivatives with respect to $x$ and $y$;
\item $I(t,d_0) := \{\bldr: D(t,\bldr) = d_0 \}$ --- the spatial level curve (isoline) of $D(t,\cdot)$ with the field level $d_0$, see Fig.~\ref{isol.fig}(a);
\item $[T,N]= [T(t,\bldr),N(t,\bldr)]$ --- the (right) Frenet frame of the spatial isoline $I[t,\gamma]$ with $\gamma := D(t,\bldr)$ at the point $\bldr$, i.e., the isoline that passes through this point at time $t$; in other words,
\begin{equation}
\label{n}
N = \frac{\nabla D(t,\bldr)}{\|\nabla D (t,\bldr)\|}
\end{equation}
and $T$ are the unit normal and tangent vectors, respectively;
\item $\varkappa(t,\bldr)$ --- the signed curvature of this spatial isoline at the point $\bldr$
($\varkappa \geq 0$ on the convexities of the domain where the field assumes greater values $\{ \bldr^\prime : D(t,\bldr^\prime) \geq \gamma \}$);
\item $\bldr_+(\Delta t| t,\bldr)$ --- the nearest (to $\bldr$) point where the $N$-axis of the above Frenet frame intersects the isoline with the same field level $\gamma$ observed at time $t=t+\Delta t$; see Fig.~\ref{fig.converge}; 
    \item $p(\Delta t| t,\bldr)$ --- the ordinate of $\bldr_+(\Delta t| t,\bldr)$, i.e., the signed normal displacement of the isoline from $t$ to $t+\Delta t$;
\item $\lambda(t,\bldr)$ --- the front velocity of the spatial isoline:
    $
\lambda (t,\bldr):=\lim_{\Delta t\to 0}\frac{p(\Delta t| t,\bldr)}{\Delta t };
$
\item $\alpha(t,\bldr)$ --- the front acceleration of the spatial isoline:
\begin{equation}
\label{alpha_def}
\alpha(t,r):= \lim_{\Delta t\to 0}\frac{\lambda[t+\Delta t, \bldr_+(\Delta t|t,\bldr)] - \lambda[t,\bldr]}{\Delta t };
\end{equation}
\item $\Delta \varphi(\Delta t|t,\bldr)$ --- the angular displacement of $T[t+\Delta t, \bldr_+(\Delta t|t,\bldr) ]$ with respect to $T[t,\bldr]$, ; see Fig.~\ref{fig.converge};
\item $\omega(t,\bldr)$ --- the angular velocity of rotation of the spatial isoline:
    $
\omega(t,\bldr):=\displaystyle{\lim_{\Delta t\to 0}}\frac{\Delta \varphi(\Delta t|t,\bldr)}{\Delta t }
$;
\item $\rho(t,\bldr)$ --- the density of isolines:
$
\rho(t,\bldr):=\displaystyle{\lim_{\Delta \gamma \to 0}}\frac{\Delta \gamma}{q(\Delta \gamma| t,\bldr)},
$
where $q(\Delta \gamma| t,\bldr)$ is the ordinate of the nearest to $\bldr$ point where the $N$-axis of the Frenet frame intersects the isoline $I(t| t,\gamma+\Delta \gamma), \gamma:= D(t,\bldr)$ with a close field value $\gamma+\Delta \gamma$; \footnote{This density characterizes the ``number'' of isolines within the unit distance from the basic one $I(t,\gamma)$, where the ``number'' is evaluated by the discrepancy in the values of $D(\cdot)$ observed on these isolines at the given time $t$.}
\item $v_\rho(t,\bldr)$ --- the (proportional) growth rate of the density $\rho$ over time:
\begin{equation}
\label{vrho.def}
v_\rho(t,\bldr):= \frac{1}{\rho(t,\bldr)}\lim_{\Delta t\to 0}\frac{\rho[t+\Delta t, \bldr_+(\Delta t|t,\bldr)]- \rho[t,\bldr]}{\Delta t};
\end{equation}
\item $ \tau_\rho(t,\bldr)$ --- the (proportional) growth rate of the density under a tangential displacement at time $t$:
\begin{equation}
\label{beta_def}
\tau_\rho (t,r):= \frac{1}{\rho(t,\bldr)} \lim_{\Delta s \to 0}\frac{\rho(t, r+ T \Delta s ) -  \rho(t,r)}{\Delta s} ;
\end{equation}
\item $ n_\rho(t,\bldr)$ --- the (proportional) growth rate of the density under a normal displacement at time $t$:
\begin{equation}
\label{njrm_def}
 n_\rho (t,r):= \frac{1}{\rho(t,\bldr)}\lim_{\Delta s \to 0}\frac{ \rho(t, r+ N \Delta s ) -  \rho(t,r)}{\Delta s} ;
\end{equation}
\item $\omega_\nabla(t,\bldr)$ --- the angular velocity of the gradient $\nabla D$ rotation at time $t$ at point $\bldr$.
\end{itemize}
The following lemma explicitly links the above quantities with the field distribution $D(\cdot)$.
\begin{lemma}
\label{lem.relation}
Whenever the field $D(\cdot)$ is twice continuously differentiable in a vicinity of $(t,\bldr)$ and $\nabla D(t,\bldr) \neq 0$, the afore-introduced quantities are well-defined and the following relations hold at $(t,\bldr)$:
 \begin{eqnarray}
 \label{speed1}
 \lambda=-\frac{D'_{t}}{\|\nabla D \|}, \quad   \rho = \|\nabla D\|, \quad  \bldr_+(dt|t,\bldr) = \bldr + \lambda N dt + \so(dt),
 \\
 \label{alpha}
 v_\rho = \frac{\spr{\nabla D^\prime_{t}+\lambda \dd{D}N}{N}}{\|\nabla D \|}, \quad
 \omega= - \frac{\spr{\nabla D^\prime_{t}+\lambda \dd{D}N}{T}}{\|\nabla D\|}, \quad
\alpha=- \frac{D^{\prime\prime}_{tt} + \lambda \left\langle \nabla D^\prime_t ; N \right\rangle}{\|\nabla D\|}  - \lambda  v_\rho ,
\\
 \label{tau}
\varkappa = - \frac{\spr{\dd{D}T}{T}}{\|\nabla D\|}, \quad
\omega_\nabla = - \frac{\spr{\nabla D^\prime_t}{T}}{\|\nabla D\|}
, \quad \tau_\rho =\frac{\spr{D^{\prime\prime}N}{T}}{\left\| \nabla D \right\|} ,
 \quad
 n_\rho =\frac{\spr{D^{\prime\prime}N}{N}}{\left\| \nabla D \right\|} .
\end{eqnarray}
\end{lemma}
\pf
We put $N:=N(t,\bldr), \gamma := D(t,\bldr), f(\Delta t,p) := D(t+\Delta t, \bldr + p N) - \gamma$, where $\Delta t, p \in \br$. Then the partial derivatives $f^\prime_p(0,0) = \spr{\nabla D(t,\bldr)}{N} = \|\nabla D(t,\bldr)\| \neq 0$ and $f^\prime_{\Delta t} (0,0) = D^\prime_t(t,\bldr)$. By the implicit function theorem \cite{KrPa02}, the equation $f(\Delta t, p) = 0$ has a unique solution $p=p(\Delta t)$ in a sufficiently small vicinity of $0$ for any sufficiently small $\Delta t$; furthermore, this solution smoothly depends on $\Delta t$ and $\frac{d p}{d(\Delta t)}(0) = - \frac{f^\prime_{\Delta t} (0,0)}{f^\prime_p(0,0)}$. This implies that $p(\Delta|t,\bldr) = p(\Delta t)$ for $\Delta t \approx 0$, the speed $\lambda$ is well defined, and the first relation from \eqref{speed1} does hold.
Similarly, the equation $\varphi (\Delta \gamma, q) := D(t, \bldr + q N) - \gamma - \Delta \gamma =0 $ has a unique solution $q=q(\Delta \gamma)$ in a sufficiently small vicinity of $0$ for any sufficiently small $\Delta \gamma$; furthermore, this solution smoothly depends on $\Delta \gamma$ and $\frac{d p}{d(\Delta \gamma)}(0) = - \frac{\varphi^\prime_{\Delta \gamma} (0,0)}{\varphi^\prime_p(0,0)} = \frac{1}{\|\nabla D(t,\bldr)\|}$.
It follows that $\rho$ is well defined and the second relation from \eqref{speed1} holds.
The third relation is immediate from the definitions of $\bldr_+(\Delta t|t,\bldr), p(\Delta t| t,\bldr)$, and $\lambda$.
\par
To proceed, we put $r_+(\Delta t):=r_+(\Delta t| t,\bldr)$. Due to \eqref{speed1},
\begin{gather}
\nonumber
\nabla D\left[t+dt,\bldr_+(dt)\right]=
\nabla D[t+dt,\bldr +\lambda N dt+\so(dt)]= \nabla D+[\nabla D^\prime_{t}+\lambda \dd{D}N]dt +\so(dt),
\\
\nonumber
\rho[t+\Delta t, \bldr_+(dt)] =
\left\| \nabla D\left[t+dt,\bldr_+(dt)\right] \right\| = \left\| \nabla D+[\nabla D^\prime_{t}+ \lambda \dd{D} N]dt +\so(dt) \right\|
\\
\nonumber
= \left\| \nabla D\right\| + \frac{\left\langle \nabla D; \nabla D^\prime_{t}+\lambda\dd{D} N \right\rangle}{\left\| \nabla D\right\|} dt + \so(dt)
\\
=
\rho + \left\langle N; \nabla D^\prime_{t}+ \lambda \dd{D} N \right\rangle dt + \so(dt) \;\mid \Rightarrow \text{the first formula in \eqref{alpha}},
\label{expan}
\\
\nonumber
N\left[t+dt,\bldr_+(dt)\right] = \frac{\nabla D\left[t+dt,\bldr_+(dt)\right]}{\left\| \nabla D\left[t+dt,\bldr_+(dt)\right] \right\|} = N + \left[ \frac{\nabla D^\prime_{t}+ \lambda \dd{D} N}{\left\| \nabla D\right\|} - \frac{\nabla D}{\left\| \nabla D\right\|^3} \left\langle \nabla D; \nabla D^\prime_{t}+ \lambda  \dd{D}N \right\rangle  \right] dt+ \so(dt)
\end{gather}
$$
=
N + \frac{1}{\left\| \nabla D\right\|}\left[\nabla D^\prime_{t}+ \lambda \dd{D} N - N \left\langle N; \nabla D^\prime_{t}+ \lambda  \dd{D}N \right\rangle  \right] dt+ \so(dt)
\overset{\text{(a)}}{=} N + \frac{\left\langle T; \nabla D^\prime_{t}+ \lambda  \dd{D}N \right\rangle} {\left\| \nabla D\right\|} T dt   + \so(dt),
$$
where (a) holds since $W = \langle W,T \rangle T + \langle W,N \rangle N \; \forall W \in \real^2$. In the Frenet frame $(T,N)$, we have
\begin{equation*}
N\left[t+dt,\bldr_+(dt)\right] = \left( \begin{array}{c}
-\sin \Delta \varphi(dt|t,\bldr)
\\
\cos \Delta \varphi(dt|t,\bldr)
\end{array}\right)
= N + \left( \begin{array}{c}
-\cos 0
\\
- \sin 0
\end{array}\right) \omega dt + \so(dt)
 =
N - \omega T  dt + \so(dt).
\end{equation*}
By equating the coefficient prefacing $T dt$ in the last two expressions, we arrive at the second formula in \eqref{alpha}. Furthermore,
\begin{multline*}
\lambda(t+dt, \bldr_+(dt)) \overset{\text{\eqref{speed1}}}{=}
-\frac{D^\prime_{t}[t+dt, \bldr_+(dt)]}{\|\nabla D(t+dt, \bldr_+(dt))\|}
\\
\overset{\text{\eqref{expan}}}{=}\lambda - \frac{D^{\prime\prime}_{tt} dt + \left\langle \nabla D^\prime_t ; \bldr_+(dt)-r \right\rangle}{\|\nabla D\|} + D^\prime_t \frac{\left\langle N; \nabla D'_{t}+\lambda D^{\prime\prime} N \right\rangle}{\|\nabla D\|^2} dt +\so(dt)
\\
\overset{\text{\eqref{speed1}}}{=} \lambda - \frac{D^{\prime\prime}_{tt} + \lambda \left\langle \nabla D^\prime_t ; N \right\rangle}{\|\nabla D\|} dt - \lambda  \frac{\left\langle N; \nabla D'_{t}+\lambda D^{\prime\prime} N \right\rangle}{\|\nabla D\|} dt +\so(dt) \overset{\text{(b)}}{=}
\lambda - \frac{D^{\prime\prime}_{tt} + \lambda \left\langle \nabla D^\prime_t ; N \right\rangle}{\|\nabla D\|} dt - \lambda  v_\rho dt +\so(dt),
\end{multline*}
where (b) follows from the first formula in \eqref{alpha}. The definition of $\alpha$ completes the proof of \eqref{alpha}. The first two equations in \eqref{tau} are well known, the third one follows from the transformation
$$
\rho(t, r+ T  ds ) = \left\| \nabla D[t, r + T ds] \right\| = \rho(t,r) + \frac{\spr{D^{\prime\prime}T}{\nabla D}}{\left\| \nabla D \right\|} ds + \so(ds)
= \rho(t,r) + \spr{D^{\prime\prime}T}{N} ds + \so(ds),
$$
the fourth equation in \eqref{tau} is established likewise. \epf
\par
It follows from Lemma~\ref{lem.relation} that $\omega = \omega_\nabla - \lambda \tau_\rho, v_\rho = - \omega_\nabla + \lambda n_\rho$.
\par
The next lemma shows how some of the introduced characteristics change under infinitesimally small shifts in time and space.
\begin{lemma}
Suppose that the field $D(\cdot)$ is twice continuously differentiable in a vicinity of $(t,\bldr)$ and $\nabla D(t,\bldr) \neq 0$.
Then the following relations hold:
\begin{gather}
\label{lrho}
\lambda (t,r+T ds) = \lambda + \omega ds + \so (ds), \quad \lambda (t,r+N ds) = \lambda - v_\rho ds + \so (ds);
\\
\label{frenet-serrat}
N[t,r+T ds] = N - \kappa T ds +\so(ds), \quad T[t,r+T ds] = T + \kappa N ds +\so(ds);
\\
\label{nnn}
N(t,r+N ds)= N + \tau_\rho T ds + \so(ds), \quad T(t,r+N ds)= T - \tau_\rho N ds + \so(ds),
\\
\label{n=omega}
N[t+dt, \bldr_+(dt)] = N - \omega T dt +\so(dt), \quad T[t+dt, \bldr_+(dt)] = T + \omega N dt +\so(dt),
\end{gather}
where $\bldr_+(\Delta t):=\bldr_+(\Delta t|t,\bldr)$.
\end{lemma}
\pf
Formulas \eqref{lrho} are justified by the following observations:
\begin{multline*}
\lambda(t,r+T ds) \overset{\text{\eqref{speed1}}}{=} -\frac{D^\prime_{t}(t,r+T ds)}{\|\nabla D(t, r+ T ds)\|} = \lambda(t,r) - \frac{\spr{\nabla D^\prime_{t}}{T}}{\|\nabla D\|} ds + D^\prime_t \frac{\spr{D^{\prime\prime}N}{T} }{\|\nabla D\|^2} ds + \so(ds)
\overset{\text{\eqref{speed1}}}{=} \lambda(t,r) - \frac{\spr{\nabla D^\prime_{t}}{T}}{\|\nabla D\|} ds
\end{multline*}
\begin{multline*}
- \lambda \frac{\spr{D^{\prime\prime}N}{T} }{\|\nabla D\|} ds + \so(ds)
=
\lambda(t,r) - \frac{\spr{\nabla D^\prime_{t}+ \lambda D^{\prime\prime}N}{T}}{\|\nabla D\|} ds  + \so(ds)
\overset{\text{\eqref{alpha}}}{=} \lambda + \omega ds + \so (ds);
\\
\lambda(t,r+N ds) = \lambda(t,r) - \frac{\spr{\nabla D^\prime_{t}+ \lambda D^{\prime\prime}N}{N}}{\|\nabla D\|} ds + \so(ds) \overset{\text{\eqref{alpha}}}{=}
\lambda - v_\rho ds + \so (ds).
\end{multline*}
Formulas \eqref{frenet-serrat} are the classic Frenet-Serrat equations. Furthermore
\begin{multline*}
N(t,r+N ds) = \frac{\nabla D(t,r+N ds)}{\|\nabla D(t,r+N ds)\|} = N + \frac{D^\pp N}{\|\nabla D\|}ds - \nabla D \frac{\spr{D^\pp N}{\nabla D}}{\|\nabla D\|^3}ds + \so(ds)
 \\
 = N + \frac{D^\pp N - N \spr{D^\pp N}{N} }{\|\nabla D\|}ds + \so(ds)
 = N + \frac{\spr{D^\pp N}{T} }{\|\nabla D\|}T ds + \so(ds) \overset{\text{\eqref{tau}}}{=} N + \tau_\rho T ds + \so(ds),
  \end{multline*}
  which yields the first formula in \eqref{nnn}. This also implies the second formula since $N = R_{\frac{\pi}{2}} T, T = - R_{\frac{\pi}{2}} N$. Formulas \eqref{n=omega} follow from the definition of $\omega$. \epf
  \par
Let us consider a simply connected domain $\mathfrak{D}$
in a vicinity of which the field $D(\cdot)$ is twice continuously differentiable and $\nabla D(t,\bldr) \neq 0$.
Then the angle $\beta$ of rotation of the vector-field $\nabla D$ along any curve lying in $\mathfrak{D}$ is uniquely determined by the ordered pair of the ends of this curve. Let $\beta_\nabla(\mathfrak{D})$ be the maximum of $|\beta|$ over all pairs in $\mathfrak{D}$.
\begin{lemma}
\label{lemrroott}
Suppose that the field $D(\cdot)$ is twice continuously differentiable and $\nabla D(t,\bldr) \neq 0$ in a vicinity of $\mathfrak{D}:=\Delta \times C$, where $T$ is an interval and $C \subset \br^2$ is a convex set. Then
$$
\beta_\nabla(\mathfrak{D}) \leq \sup_{p \in \mathfrak{D}} |\omega_\nabla(p)| \times |T| + \sup_{p \in \mathfrak{D}} \sqrt{\varkappa^2+\tau_\rho^2} \times \text{\bf diam} C,
$$
where $|T|$ is the length of the interval, $\text{\bf diam}$ is the diameter of the set, i.e., the supremum of distances between two its elements,
and $\omega_\nabla, \varkappa$, and $\tau_\rho$ are the field parameters introduced at the beginning of this section.
\end{lemma}
\pf
Let $p_i=(t_i, \bldr_i) \in \mathfrak{D}:=\Delta \times C, i=0,1$ and $t_0 \leq t_1$ for the definiteness. In $\mathfrak{D}$, we introduce two parametric curves $\zeta(s):= [s, \bldr_0], s \in [t_0,t_1]$ and $\xi(s):= [t_1, (1-s)\bldr_0 + s \bldr_1], s \in [0,1]$, connecting the point $[t_0,\bldr_0]$ with $[t_0,\bldr_1]$ and $[t_0,\bldr_1]$ with $[t_1,\bldr_1]$, respectively. The concatenation of these curves connects $p_1$ with $p_2$ and lies in $\mathfrak{D}$. So the angle of the gradient rotation when going from $p_1$ to $p_2$ inside $\mathfrak{D}$ is the angle of rotation $\sphericalangle \beta_\zeta $ when going along $\zeta(\cdot)$ plus the angle of rotation $\sphericalangle \beta_\xi $ when going along $\xi(\cdot)$. Let $\varphi[t,\bldr]$ stands for the orientation angle of $\nabla D [t,\bldr]$ in the absolute coordinate frame. Then
\begin{eqnarray*}
\frac{d}{ds} \varphi[\zeta(s)] = - \frac{\spr{\nabla D^\prime_t [\zeta(s)]}{T[\zeta(s)]}}{\|\nabla D[\zeta(s)]\|} \overset{\text{\eqref{tau}}}{=\!=} \omega_\nabla[\zeta(s)] \Mapsto |\sphericalangle \beta_\zeta | = \left| \int_{t_0}^{t_1} \omega_\nabla[\zeta(s)] \;ds \right| \leq \sup_{p \in \mathfrak{D}} |\omega_\nabla(p)| (t_1-t_0);
\\
\frac{d}{ds} \varphi[\xi(s)] = - \frac{\spr{D^\pp [\zeta(s)][\bldr_1-\bldr_0]}{T[\zeta(s)]}}{\|\nabla D[\zeta(s)]\|}
= - \frac{\spr{D^\pp [\zeta(s)]T[\zeta(s)]}{T[\zeta(s)]}}{\|\nabla D[\zeta(s)]\|} \spr{\bldr_1-\bldr_0}{T[\zeta(s)]} \\ - \frac{\spr{D^\pp [\zeta(s)]T[\zeta(s)]}{N[\zeta(s)]}}{\|\nabla D[\zeta(s)]\|} \spr{\bldr_1-\bldr_0}{N[\zeta(s)]} \overset{\text{\eqref{tau}}}{=\!=}
 +\varkappa \spr{\bldr_1-\bldr_0}{T[\zeta(s)]} - \tau_\rho \spr{\bldr_1-\bldr_0}{N[\zeta(s)]} \Mapsto
\\
\left| \frac{d}{ds} \varphi[\xi(s)] \right| \leq \|\bldr_1-\bldr_0\| \sqrt{\varkappa^2+\tau_\rho^2} \Mapsto |\sphericalangle \beta_\xi | \leq
\int_0^1 \left| \frac{d}{ds} \varphi[\xi(s)] \right| \; ds \leq \|\bldr_1-\bldr_0\| \sup_{p \in \mathfrak{D}} \sqrt{\varkappa^2+\tau_\rho^2}.
\end{eqnarray*}
To complete the proof, it remains to take the supremum over $p_0,p_1 \in \mathfrak{D}$. \epf
\section{Time Derivatives of the Robot's Field Readings}
\label{robot.deriv}
\setcounter{equation}{0}
\begin{theorem}
Let the robot \eqref{1aa} move at the speed $v \equiv \ov{v}$ in a domain where the field $D(\cdot)$ is twice continuously differentiable and $\nabla D(t,\bldr) \neq 0$. Then the following relations hold for the field reading $d(t):=D[t,\bldr(t)]$:
\begin{gather}
\label{for.d.d}
\dot{d} = \rho \left[ \ov{v} \spr{N}{\vec{e}} - \lambda\right],
\\
\label{for.d.dd}
\ddot{d}(t) =  \rho \left[(\dot{\theta} - 2\omega  - \varkappa v_T + 2v_\Delta \tau_\rho ) v_T   +2 v_\rho v_\Delta  - \alpha + v_\Delta^2 n_\rho
 \right].
\end{gather}
Here $\rho, \lambda, v_\rho, n_\rho, \tau_\rho, \alpha, \omega, \varkappa$ are the field parameters introduced at the beginning of Section~{\rm \ref{subsec.not}}, whereas
\begin{equation}
\label{def.vind}
v_{\Delta} := \frac{\dot{d}}{\rho} \overset{\text{\eqref{for.d.d}}}{=\!=} \ov{v} \spr{N}{\vec{e}} - \lambda =  \spr{N}{\vec{v}} - \lambda, \qquad
v_T := \spr{\vec{v}}{T} = \sgn \sigma \sqrt{\ov{v}^2 -\left[ \lambda  + v_\Delta \right]^2}, \quad \text{\rm where} \; \sigma := \spr{\vec{e}}{T}
\end{equation}
and the Frenet frame $[T,N]$ is built at the robot's current location.
\end{theorem}
\pf Formula \eqref{for.d.d} is true since
$$
\dot{d} \overset{\text{\eqref{1aa}}}{=\!=} \ov{v} \spr{\nabla D}{\vec{e}} +  D^\prime_t = \|\nabla D\| \left[ \ov{v} \spr{\frac{\nabla D}{\|\nabla D\|}}{\vec{e}} +  \frac{D^\prime_t}{|\nabla D\|} \right] \overset{\text{(a)}}{=} \rho \left[ \ov{v} \spr{N}{\vec{e}} - \lambda\right],
$$
where (a) holds thanks to \eqref{n} and the first two formulas in \eqref{speed1}.
\par
To establish \eqref{for.d.dd}, we first note that due to \eqref{def.vind}, $\spr{N}{\vec{v}} = v_{\Delta} + \lambda$ and
\begin{equation}
\label{form.vecv}
\vec{v} = \spr{T}{\vec{v}} T + \spr{N}{\vec{v}} N = v_T T + (v_{\Delta} + \lambda) N, \qquad \vec{e} = \frac{\vec{v}}{\ov{v}} =\frac{v_T}{\ov{v}} T + \frac{v_{\Delta} + \lambda}{\ov{v}} N.
\end{equation}
This justifies (a) in the following stream of transformations:
\begin{multline}
 r(t +dt) \overset{\text{\eqref{1}}}{=\!=} r(t) + \vec{v} dt +\so(dt) \overset{\text{(a)}}{=} \underbrace{r(t) + \lambda N \, dt}_{= \bldr_+(dt)+
 \text{\tiny $\mathcal{O}$}(dt) \text{by \eqref{speed1}}}
 +  [ v_\Delta N +  v_T  T ] dt +\so(dt)
 \\
 = \bldr_+(dt)+ [ v_\Delta N +  v_T  T ] dt +\so(dt) .
 \label{rplus}
 \end{multline}
 We proceed by using \eqref{for.d.d} and dropping the argument $[t,\bldr(t)]$ everywhere:
 \begin{multline}
 \dot{d}(t+dt) - \dot{d}(t) =
  \rho[t+dt,\bldr(t+dt)] \Big( - \lambda [t+dt,\bldr(t+dt)]
  + \ov{v} \left\langle N[t+dt,\bldr(t+dt)], \vec{e}(t+dt) \right\rangle \Big)
\\
 - \rho \left( - \lambda + \ov{v} \left\langle N, \vec{e} \right\rangle \right)
 \\
 = \big\{ \rho[t+dt,\bldr(t+dt)] - \rho\big\}\underbrace{\left( - \lambda [t+dt,\bldr(t+dt)]
  + \ov{v} \left\langle N[t+dt,\bldr(t+dt)], \vec{e}(t+dt) \right\rangle\right)}_{=  - \lambda + \ov{v} \left\langle N, \vec{e} \right\rangle + \eta(dt),\; \text{where}\; \eta(z) \to 0 \; \text{as}\; z \to 0}
 \\
 - \rho \left\{ \lambda [t+dt,\bldr(t+dt)] - \lambda \right\}
 \\
 + \ov{v} \rho  \big\langle N[t+dt,\bldr(t+dt)] - N, \vec{e} \big\rangle
 + \ov{v} \rho \big\langle N, \vec{e}(t+dt) - \vec{e} \big\rangle 
\\
 = \underbrace{\big\{ \rho[t+dt,\bldr(t+dt)] - \rho\big\}}_{=:a}\underbrace{( - \lambda + \ov{v} \left\langle N, \vec{e} \right\rangle)}_{= v_\Delta\; \text{by \eqref{def.vind}}}
 - \rho \underbrace{\left\{ \lambda [t+dt,\bldr(t+dt)] - \lambda \right\}}_{=:b}
 \\
 + \ov{v} \rho  \big\langle \underbrace{N[t+dt,\bldr(t+dt)] - N}_{=:\vec{c}_1}, \vec{e} \big\rangle
 + \ov{v} \rho \big\langle N, \underbrace{\vec{e}(t+dt) - \vec{e}}_{=:\vec{c}_2} \big\rangle + \so(dt).
 \label{first,expan}
 \end{multline}
 Now we analyze $a,b,\vec{c}_1$, and $\vec{c}_2$ separately:
\begin{multline*}
 a \overset{\text{\eqref{rplus}}}{=\!=}  \rho[t+dt,\bldr_+(dt) + v_\Delta N dt + v_T T dt + \so(dt)] - \rho
 \\
 =  \left\{\rho[t+dt,\underbracket{\bldr_+(dt) + v_\Delta N dt + \so(dt)} + v_T T dt] - \rho[t+dt,\underbracket{\bldr_+(dt) + v_\Delta N dt + \so(dt)}] \right\}
 \\
 + \rho[t+dt,\bldr_+(dt) + v_\Delta N dt + \so(dt)] - \rho
 \\
  \overset{\text{\eqref{beta_def}}}{=\!=} \left\{ \rho v_T \tau_\rho dt + \so(dt) \right\} +  \rho[t+dt,\bldr_+(dt) + v_\Delta N dt + \so(dt)] - \rho
  \\
=   \rho v_T \tau_\rho dt  + \left\{ \rho[t+dt,\underbracket{\bldr_+(dt) + \so(dt)} + v_\Delta N dt ] - \rho[t+dt,\underbracket{\bldr_+(dt) + \so(dt)}]\right\}
\\
+ \rho[t+dt,\bldr_+(dt) + \so(dt)] - \rho + \so(dt)
\\
\overset{\text{\eqref{njrm_def}}}{=\!=} \rho v_T \tau_\rho dt  + \left\{ \rho v_{\Delta} n_\rho dt +  \so(dt) \right\} + \rho[t+dt,\bldr_+(dt) + \so(dt)] - \rho + \so(dt)
\\
= \rho [ v_T \tau_\rho +v_{\Delta} n_\rho ]dt  + \rho[t+dt,\bldr_+(dt)] - \rho + \so(dt)
\\
 \overset{\text{\eqref{vrho.def}}}{=\!=}
 \rho \cdot [v_\rho + v_\Delta n_\rho + v_T \tau_\rho ]dt + \so(dt).
\end{multline*}
Similarly,
\begin{multline*}
\nonumber
b = \lambda [t+dt,\bldr(t+dt)] - \lambda \overset{\text{\eqref{rplus}}}{=\!=} \lambda \left[t+dt,\bldr_+(dt) + v_\Delta N dt + v_T T dt + \so(dt) \right] - \lambda
\\
= \left\{\lambda[t+dt,\underbracket{\bldr_+(dt) + v_\Delta N dt + \so(dt)} + v_T T dt] - \lambda[t+dt,\underbracket{\bldr_+(dt) + v_\Delta N dt + \so(dt)}] \right\}
 \\
 + \lambda[t+dt,\bldr_+(dt) + v_\Delta N dt + \so(dt)] - \lambda
 \\
  \overset{\text{\eqref{lrho}}}{=\!=} \left\{ v_T \omega dt + \so(dt) \right\} +  \lambda[t+dt,\bldr_+(dt) + v_\Delta N dt + \so(dt)] - \lambda
    \\
=  v_T \omega dt  + \left\{ \lambda[t+dt,\underbracket{\bldr_+(dt) + \so(dt)} + v_\Delta N dt ] - \lambda[t+dt,\underbracket{\bldr_+(dt) + \so(dt)}]\right\}
\\
+ \lambda[t+dt,\bldr_+(dt) + \so(dt)] - \lambda + \so(dt)
\\
\overset{\text{\eqref{lrho}}}{=\!=} v_T \omega dt  + \left\{ - v_{\Delta} v_\rho dt +  \so(dt) \right\} + \lambda[t+dt,\bldr_+(dt) + \so(dt)] - \lambda + \so(dt)
\\
= [ v_T \omega -v_{\Delta} v_\rho ]dt  + \lambda[t+dt,\bldr_+(dt)] - \lambda + \so(dt)
\\
\overset{\text{\eqref{alpha_def}}}{=\!=} [\alpha + \omega v_T - v_\rho v_\Delta] dt + \so(dt);
\end{multline*}
\begin{multline*}
\vec{c}_1= N[t+dt,\bldr(t+dt)] - N \overset{\text{\eqref{rplus}}}{=\!=} N \left[t+dt,\bldr_+(dt) + v_\Delta N dt + v_T T dt + \so(dt) \right] - N
\\
= \left\{N[t+dt,\underbracket{\bldr_+(dt) + v_\Delta N dt + \so(dt)} + v_T T dt] - N[t+dt,\underbracket{\bldr_+(dt) + v_\Delta N dt + \so(dt)}] \right\}
 \\
 + N[t+dt,\bldr_+(dt) + v_\Delta N dt + \so(dt)] - N
 \\
  \overset{\text{\eqref{frenet-serrat}}}{=\!=} \left\{ - v_T \varkappa T dt + \so(dt) \right\} +  N[t+dt,\bldr_+(dt) + v_\Delta N dt + \so(dt)] - N
      \\
=  - v_T \varkappa T dt  + \left\{ N[t+dt,\underbracket{\bldr_+(dt) + \so(dt)} + v_\Delta N dt ] - N[t+dt,\underbracket{\bldr_+(dt) + \so(dt)}]\right\}
\\
+ N[t+dt,\bldr_+(dt) + \so(dt)] - N + \so(dt)
\\
\overset{\text{\eqref{nnn}}}{=\!=} - v_T \varkappa dt  + \left\{ v_{\Delta} \tau_\rho T dt +  \so(dt) \right\} + N [t+dt,\bldr_+(dt) + \so(dt)] - N + \so(dt)
\\
= [ v_{\Delta} \tau_\rho - v_T \varkappa]T dt  + N[t+dt,\bldr_+(dt)] - N + \so(dt)
\\
\overset{\text{\eqref{n=omega}}}{=\!=} [- \omega + \tau_\rho v_\Delta - \varkappa v_T] T dt + \so(dt).
\end{multline*}
Finally,
\begin{multline*}
\vec{c}_2 = \vec{e}(t+dt) - \vec{e} \overset{\text{\eqref{1aa}}}{=\!=} \dot{\theta} R_{\frac{\pi}{2}} \vec{e} dt +\so(dt) =
\frac{\dot{\theta}}{\ov{v}} R_{\frac{\pi}{2}} \vec{v} dt +\so(dt) \overset{\text{\eqref{form.vecv}}}{=\!=}
\frac{\dot{\theta}}{\ov{v}} R_{\frac{\pi}{2}} \left[ v_T T + (v_{\Delta} + \lambda) N \right] dt +\so(dt)
\\
= \frac{\dot{\theta}}{\ov{v}} \left[ v_T R_{\frac{\pi}{2}}  T + (v_{\Delta} + \lambda) R_{\frac{\pi}{2}}  N \right] dt +\so(dt)
= \frac{\dot{\theta}}{\ov{v}} \left[- \left( \lambda   + v_\Delta \right) T + v_T  N \right] dt + \so(dt).
\end{multline*}
Now we put the obtained expressions for $a,b,\vec{c}_1, \vec{c}_2$ into \eqref{first,expan} and note that the left hand side of \eqref{first,expan}
equals $\dot{d}(t+dt) - \dot{d}(t) = \ddot{d}(t) dt + \so(dt)$. As a result, we see that
\begin{multline*}
\ddot{d}(t) = \underbrace{\rho \cdot [v_\rho + v_\Delta n_\rho + v_T \tau_\rho ]}_{\sim a}v_\Delta
 - \rho \underbrace{[\alpha + \omega v_T - v_\rho v_\Delta]}_{\sim b}
 \\
 + \ov{v} \rho  \big\langle \underbrace{[- \omega + \tau_\rho v_\Delta - \varkappa v_T]T}_{\sim \vec{c}_1}, \vec{e} \big\rangle
 + \ov{v} \rho \Big\langle N, \underbrace{\frac{\dot{\theta}}{\ov{v}} \left[- \left( \lambda   + v_\Delta \right) T + v_T  N \right]}_{\sim \vec{c}_2} \Big\rangle
 \\
 \overset{\text{\eqref{form.vecv}}}{=\!=}
 \rho \left[ (v_\rho + v_\Delta n_\rho + v_T \tau_\rho )v_\Delta
 - (\alpha + \omega v_T - v_\rho v_\Delta) \right]
 \\
 + \ov{v} \rho  \left\langle [- \omega + \tau_\rho v_\Delta - \varkappa v_T]T; \frac{v_T}{\ov{v}} T + \frac{v_{\Delta} + \lambda}{\ov{v}} N \right\rangle
 + \rho \dot{\theta}\Big\langle N ; - \left( \lambda   + v_\Delta \right) T + v_T  N \Big\rangle
 \\
=  \rho \left[ (v_\rho + v_\Delta n_\rho + v_T \tau_\rho )v_\Delta
 - (\alpha + \omega v_T - v_\rho v_\Delta) \right]
 + v_T \rho  (- \omega + \tau_\rho v_\Delta - \varkappa v_T)
 + v_T \rho \dot{\theta}
 \\
 =
 \rho \left[ (v_\rho + v_\Delta n_\rho + v_T \tau_\rho )v_\Delta
 - (\alpha + \omega v_T - v_\rho v_\Delta) + v_T  (- \omega + \tau_\rho v_\Delta - \varkappa v_T) + v_T \dot{\theta} \right]
 \\
 =
 \rho \left[(\dot{\theta} - 2\omega  - \varkappa v_T + 2v_\Delta \tau_\rho ) v_T   +2 v_\rho v_\Delta  - \alpha + v_\Delta^2 n_\rho
 \right],
\end{multline*}
which completes the proof of \eqref{for.d.dd}.
\epf
\section{Deviation of a Perpetually Rotating Robot from the Initial State}
\label{init.sec.state}
\setcounter{equation}{0}
\begin{lemma}
\label{lem.diberg}
 Let the robot move at the maximal speed $v \equiv \ov{v}$ and so that the vector $\vec{e}(\theta)$ rotates in a constant direction with the angular speed $ |\dot{\theta}| \geq \omega_\theta$. Then its deviation from the initial location obeys the inequality
\begin{equation}
\label{ineq.dist}
\|\bldr(t) - \bldr(0)\|
\leq q(\varphi):= \frac{2\ov{v}}{\omega_\theta} \left\lfloor \frac{\varphi}{2 \pi} \right\rfloor + \frac{\ov{v}}{\omega_\theta} \left[ 1 - \cos \min\{ \Lbag \varphi \Rbag; \pi \}\right] \leq \frac{2\ov{v}}{\omega_\theta} \left\lceil \frac{\varphi}{2 \pi} \right\rceil,
\end{equation}
where $\varphi := |\theta(t)-\theta(0)|,  \Lbag \varphi \Rbag := \varphi - 2 \pi \left\lfloor \frac{\varphi}{2 \pi} \right\rfloor$, whereas $\lfloor z \rfloor$ and $\lceil z \rceil$ are the integer floor and ceiling of a real number $z$, respectively.
\end{lemma}
\pf
Without any loss of generality, it can be assumed that $\bldr(0)=0, \theta(0) =0$.
We focus on the case, where $\vec{e}$ rotates clockwise, the case of counter clockwise rotation is considered likewise.
Then
$
 \dot{\theta} \leq  - \omega_\theta <0$ and so
$s:= -\theta$ can be taken as a new independent variable:
$$
\frac{d \bldr }{d s} = u \vec{e}(-s), \quad u:= - \frac{\ov{v}}{\dot{\theta}} \in \left[ 0 , \frac{\ov{v}}{\omega_\theta} \right].
$$
The squared deviation does not exceed the maximal value of $\mathcal{I}$ in the following optimization problem:
\begin{equation*}
\mathcal{I}:= \|\bldr(\varphi)\|^2 \to \max \quad \text{subject to} \quad \frac{d \bldr }{d s} = u \vec{e}(-s) \; s \in [0,\varphi], \quad \bldr(0)=0, \quad u (s) \in \left[ 0 , \ov{v}/\omega_\theta \right].
\end{equation*}
Its solution $\bldr^0(\cdot), u^0(\cdot)$ exists and obeys the Pontryagin's maximum principle \cite{Vin00}: There exists a smooth function $\psi (s) \in \br^2$ of variable $s \in [0,\varphi]$ such that
$$
u^0(s) = \text{arg\,max}_{u \in [0, \ov{v}/\omega_\theta]} u \psi^\trs \vec{e}(-s)
= \begin{cases}
 \ov{v}/\omega_\theta & \text{if}\; \psi^\trs \vec{e}(-s) >0
 \\
 0 & \text{if}\; \psi^\trs \vec{e}(-s) < 0
  \\
 \text{unclear} & \text{if}\; \psi^\trs \vec{e}(-s) = 0
\end{cases},
$$
where
$$
\frac{d \psi }{ds} = - \frac{\partial}{\partial \bldr} u^0 \psi^\trs \vec{e}(-s) =0 \Rightarrow \psi = \text{const},
\qquad
\psi = \psi(\varphi) = 2 \bldr^0(\varphi).
$$
If $\bldr^0(\varphi) = 0$, \eqref{ineq.dist} is evident.
Let $\bldr^0(\varphi)\neq 0$. Then $\psi \neq 0$ and
so as $s$ progresses, $u^0(s)$ interchanges the values $0$ and $\ov{v}/\omega_\theta$, each taken on an interval of length $\pi$ possibly except for  extreme intervals whose lengths do not exceed $\pi$. Inequality \eqref{ineq.dist} results from computation of $\|\bldr(\varphi)\|$ for such controls $u(\cdot)$, along with picking the maximum among these results.
\epf

\bibliographystyle{plain}
  \bibliography{Hamidref}
 \end{document}